\documentclass[runningheads]{lncse}

\usepackage{amsmath}
\usepackage{amsfonts, amscd, stmaryrd, wasysym, enumerate, amssymb, color, graphicx, xspace, verbatim, mathtools, paralist}
\usepackage{epsfig}
\usepackage{shuffle}
\usepackage{difftrees}
\usepackage{graphics} 

\newcommand{\tr}{\triangleright}
\newcommand{\id}{{\mathrm{id}}}
\newcommand{\OT}{\operatorname{OT}}

\def\g{{\mathfrak g}}

\def\so{{\mathfrak s\mathfrak o}}
\def\RR{{\mathbb R}}
\def\postLie{\operatorname{postLie}}
\def\F{\mathcal F}
\def\one{{\mathbb I}}
\newcommand{\OF}{\text{OF}}
\def\glb{\g_{\operatorname{LB}}}
\def\Glb{G_{\operatorname{LB}}}

\begin{document}


\title{What is a post-Lie algebra and why is it useful in geometric integration}

\titlerunning{Post-Lie algebra in geometric integration}

\author{Charles Curry\inst{1} \and Kurusch Ebrahimi-Fard\inst{1} \and Hans Munthe-Kaas\inst{2}}

\authorrunning{Curry et al.}   

\institute{
Norwegian University of Science and Technology (NTNU), Institutt for matematiske fag, 
N-7491 Trondheim, Norway {\tt charles.curry@ntnu.no, kurusch.ebrahimi-fard@ntnu.no}
\and 
University of Bergen, Department of Mathematics, N-5020 Bergen, Norway {\tt hans.munthe-kaas@uib.no}
}

\maketitle


\begin{abstract}
We explain the notion of a post-Lie algebra and outline its role in the theory of Lie group integrators. 
\end{abstract}


\section{Introduction}
\label{sect:intro}

In recent years classical numerical integration methods have been extended beyond applications in Euclidean space onto manifolds. In particular, the theory of Lie group methods \cite{Iserles00} has been developed rapidly. In this respect Butcher's $B$-series \cite{HWL} have been generalized to Lie--Butcher series \cite{MK1,MK2}. Brouder's work \cite{Brouder} initiated the unfolding of rich algebro-geometric aspects of the former, where Hopf and pre-Lie algebras on non-planar rooted trees play a central role \cite{CHV,MMMKV}. Lie--Butcher series underwent similar developments replacing non-planar trees by planar ones \cite{LMK,MKW}. Correspondingly, pre-Lie algebras are to $B$-series what post-Lie algebras are to Lie-Butcher series \cite{EFLMK,FMK}.   

In this note we explore the notion of a post-Lie algebra and outline its importance to integration methods.


\section{Post-Lie algebra and examples}
\label{sect:postlie}

We begin by giving the definition of a post-Lie algebra followed by a proposition describing the central result. The three subsequent examples illustrate the value of such algebras, in particular to Lie group integration methods. 

\begin{definition} \label{def:postLie} 
A \emph{post-Lie algebra} $(\mathfrak{g}, [\cdot,\cdot], \triangleright)$ consists of a Lie algebra $(\mathfrak g, [\cdot,\cdot])$ and a binary product $\triangleright  : {\mathfrak g} \otimes {\mathfrak g} \rightarrow \mathfrak g$ such that, for all elements $x,y,z \in \mathfrak g$ the following relations hold
\begin{align}
\label{postLie1}
	x \triangleright [y,z] &= [x\triangleright y , z] + [y , x \triangleright z],\\
\label{postLie2}
	[x,y] \triangleright z &= {\rm{a}}_{\triangleright  }(x,y,z) - {\rm{a}}_{\triangleright}(y,x,z),
\end{align}
where the associator ${\rm{a}}_{\triangleright}(x,y,z):=x \triangleright (y \triangleright z) - (x \triangleright y) \triangleright z$.  
\end{definition}

Post-Lie algebras first appear in the work of Vallette \cite{Vallette} and were independently described in \cite{MKW}. Comparing these references the reader will quickly see how different they are in terms of aim and scope, which hints at the broad mathematical importance of this structure.

\begin{proposition}\label{prop:Liebracket}
Let $(\mathfrak g, [\cdot, \cdot], \triangleright )$ be a post-Lie algebra. For $x, y \in \mathfrak g$ the bracket
\begin{equation}
\label{postLie3}
	\llbracket x,y \rrbracket := x \tr y - y \tr x + [x,y]
\end{equation}
satisfies the Jacobi identity. The resulting Lie algebra is denoted $(\overline{\mathfrak g},\llbracket \cdot,\cdot \rrbracket)$. 
\end{proposition}

\begin{corollary}\label{cor:preLie}
A post-Lie algebra with an abelian Lie algebra $(\mathfrak g, [\cdot, \cdot]=0, \triangleright )$ reduces to a left pre-Lie algebra, i.e., for all elements $x,y,z \in \mathfrak g$ we have
\begin{align}
\label{preLie}
	{\rm{a}}_{\triangleright  }(x,y,z) = {\rm{a}}_{\triangleright}(y,x,z).
\end{align}
\end{corollary}

\begin{example}\label{exm:connection}
Let $\mathcal{X}(M)$ be the space of vector fields on a manifold $M$, equipped with a linear connection. The covariant derivative $\nabla_X Y$ of $Y$ in the direction of $X$ defines an $\mathbb{R}$-linear, non-associative binary product $X\triangleright Y$ on $\mathcal{X}(M)$. The torsion $T$, a skew-symmetric tensor field of type $(1,2)$, is {defined by}
\begin{align}
\label{torsion}
	T(X,Y) := X\triangleright Y - Y\triangleright X - \llbracket X,Y \rrbracket,
\end{align}
where the bracket $\llbracket \cdot,\cdot \rrbracket$ on the right is the Jacobi bracket of vector fields. The torsion admits a covariant differential $\nabla T$, a tensor field of type $(1,3)$. Recall the definition of the curvature tensor $R$, a tensor field of type $(1,3)$ given by
\[
	R(X,Y)Z = X\triangleright(Y\triangleright Z) - Y\triangleright(X\triangleright Z) 
	- \llbracket X,Y \rrbracket \triangleright Z.
\]
In the case that the connection is flat and has constant torsion, i.e., $R=0=\nabla T$, we have that $(\mathcal{X}(M),-T(\cdot,\cdot),\triangleright)$ defines a post-Lie algebra. Indeed, the first Bianchi identity shows that $-T(\cdot,\cdot)$ obeys the Jacobi identity; as $T$ is skew-symmetric it therefore defines a Lie bracket. Moreover, flatness is equivalent to $(\ref{postLie2})$ as can be seen by inserting \eqref{torsion} into the statement $R=0$, whilst $(\ref{postLie1})$ follows from the definition of the covariant differential of $T$:
\[
	0 = \nabla T (Y,Z;X) = X \triangleright T(Y,Z) - T(Y,X \triangleright Z) - T(X \triangleright Y,Z).
\]
The formalism of post-Lie algebras assists greatly in understanding the interplay between covariant derivatives and integral curves of vector fields, which is central to the study of numerical analysis on manifolds.
\end{example}


\begin{example}\label{exm:rootedtrees} We now consider planar rooted trees with left grafting. Recall that a rooted tree is made out of vertices and non-intersecting oriented edges, such that all but one vertex have exactly one outgoing line and an arbitrary number of incoming lines. The root is the only vertex with no outgoing line and is drawn on bottom of the tree, whereas the leaves are the only vertices without any incoming lines. A planar rooted tree is a rooted tree with an embedding in the plane, that is, the order of the branches is fixed. We denote the set of planar rooted trees by $\OT$. 
\[
	\OT = \Big\{\begin{array}{c}
		\scalebox{0.6}{\ab}, \scalebox{0.6}{\aabb},\scalebox{0.6}{\aaabbb}, 
		\scalebox{0.6}{\aababb}, \scalebox{0.6}{\aaaabbbb},
		\scalebox{0.6}{\aaababbb},\scalebox{0.6}{\aaabbabb},
		\scalebox{0.6}{\aabaabbb}, \scalebox{0.6}{\aabababb},\ldots
 			\end{array}
			\Big\}.
\]
The left grafting of two trees $\tau_1 \triangleright \tau_2$ is the sum of all trees resulting from attaching the root of $\tau_1$ via a new edge successively to all the nodes of the tree $\tau_2$ from the left.
\begin{align}
  	\scalebox{0.6}{\AabB}\,\triangleright \scalebox{0.6}{\aababb}
	=\scalebox{0.6}{\aAabBababb} + \scalebox{0.6}{\aaAabBbabb} + \scalebox{0.6}{\aabaAabBbb}.
\end{align}
Left grafting means that the tree $\tau_1$, when grafted to a vertex of $\tau_2$ becomes the leftmost branch of this vertex. We consider now the free Lie algebra $\mathcal{L}(\OT)$ generated by planar rooted trees. In \cite{LMK} is was shown that $\mathcal{L}(\OT)$ together with left grafting defines a post-Lie algebra. In fact, it is the free post-Lie algebra $\mathrm{PostLie}(\scalebox{0.6}{\ab})$ on one generator \cite{LMK}. 

Ignoring planarity, that is, considering non-planar rooted trees, turns left grafting into grafting, which is a pre-Lie product on rooted tree satisfying \eqref{preLie} \cite{Manchon}. The space spanned by non-planar rooted trees together with grafting defines the free pre-Lie algebra $\mathrm{PreLie}(\scalebox{0.6}{\ab})$ on one generator \cite{CL}.
\end{example}

\begin{example}\label{exm:rmatrices} 
Another rather different example of post-Lie algebra comes from projections on the algebra $\mathcal{M}_n(\mathbb{K})$ of $n \times n$ matrices with entries in the base field $\mathbb{K}$. More precisely, we consider linear projections involved in classical matrix factorization schemes, such as $LU$, $QR$ and Cholesky \cite{CN,DLT}. Let $\pi^\ast_+$ be such a projection on $\mathcal{M}_n(\mathbb{K})$, where $\ast=LU$, $QR$, $Ch$. It turns out that both $\pi^\ast_+$ and $\pi^\ast_-:=\id -\pi^\ast_+$ satisfy the Lie algebra identity    
$$
	[\pi^\ast_\pm M,\pi^\ast_\pm N]+\pi^\ast_\pm [M,N]
	=\pi^\ast_\pm([\pi^\ast_\pm M,N] + [M,\pi^\ast_\pm N]),
$$      
for all $M,N \in \mathcal{M}_n(\mathbb{K})$. In \cite{BGN} it was shown that $M \tr N := -[\pi^\ast_- M,N] $ defines a post-Lie algebra with respect to the Lie algebra defined on $\mathcal{M}_n(\mathbb{K})$. Corollary \ref{cor:preLie} is more subtle in this context as it reflects upon the difference between classical and modified classical Yang--Baxter equation \cite{DLT,EFMMK,EFM}.        
\end{example}


\section{Post-Lie algebras and Lie group integration} 
\label{sect:integration}

We now consider post-Lie algebras as they appear in numerical Lie group integration. Recall the standard formulation of Lie group integrators~\cite{Iserles00}, where differential equations on a homogeneous space $M$ are formulated using a left action $\cdot\colon G\times M\rightarrow M$ of a Lie group $G$ of isometries on $M$, with Lie algebra $\mathfrak{g}$. An infinitesimal action $\cdot\colon \g\times M\rightarrow TM$ arises from differentiation,
\[
	V\cdot p := \left.\frac{\partial}{\partial t}\right|_{t=0}\exp(tV)\cdot p .
\]
In this setting any ordinary differential equation on $M$ can be written as
\begin{equation}
	y'(t) = f(y(t))\cdot y(t),\label{eq:ode}
\end{equation}
where $f\colon M\rightarrow \g$. For instance ODEs on the 2-sphere $S^2\simeq SO(3)/SO(2)$ can be expressed using the infinitesimal action of $\so(3)$. Embedding $S^2\subset \mathbb{R}^3$ realizes the action and infinitesimal action as matrix-vector multiplications, where $SO(3)$ is the space of orthogonal matrices, and $\so(3)$ the skew-symmetric matrices.

\smallskip

To obtain a description of the solution of (\ref{eq:ode}), we begin by giving a post-Lie algebra structure to $\g^M$, the set of (smooth) functions from $M$ to $\g$. For $f,g\in \g^M$, we let $[f,g](p) := [f(p),g(p)]_\g$ and 
$$
	(f\tr g)(p) := \left.\frac{\partial}{\partial t}\right|_{t=0} g(\exp(tf(p))\cdot p).
$$
The flow map of (\ref{eq:ode}) admits a Lie series expansion, where the terms are differential operators of arbitrary order, which live in the enveloping algebra of the Lie algebra generated by the infinitesimal action of $f$. Recall that for a Lie algebra $(\g,[\cdot,\cdot])$, the enveloping algebra is an associative algebra $(U(\g),\cdot)$ such that $\g\subset U(\g)$ and $[x,y] = x\cdot y-y\cdot x$ in $U(\g)$. As a Lie algebra $\g$ with product $\tr$, the enveloping algebra of $(\mathfrak{\g}, [\cdot,\cdot], \triangleright)$ is $U(\g)$ together with an extension of $\tr$ onto $U(\g)$ defined such that for all $x \in \g$ and $y,z\in U(\g)$
\begin{align*}
	x\tr (y\cdot z) &= (x\tr y)\cdot z + y\cdot (x\tr z)\\
	(x\cdot y)\tr z &= a_\tr(x,y,z).
\end{align*}
Many of these operations are readily computable in practice. Recall that $\g$ has a second Lie algebra structure $\bar{\g}$ associated to the bracket $\llbracket\cdot,\cdot\rrbracket$, reflecting the Jacobi bracket of the vector fields on $M$ generated by the infinitesimal action of $\g$. As a vector space, its enveloping algebra $\left(U(\bar{\g}),\ast\right)$ is isomorphic to $U(\g)$. The Lie series solution of (\ref{eq:ode}) is essentially the exponential in $U(\bar{\g})$, which in contrast to the operations of $U(\g)$ is in general difficult to compute. We are lead to the following:\newline
\newline
\noindent {\it{Basic aim:}} The fundamental problem of numerical Lie group integration is the approximation of the exponential $\exp^{\ast}$ in $\left(U(\bar{\g}),\ast\right)$ in terms of the operations of $\left(U(\g),\cdot,\tr\right)$, where $(\mathfrak{g}, [\cdot,\cdot], \triangleright)$ is the free Post-Lie algebra over a single generator.
\begin{remark}
One may wonder why we use post-Lie algebras, which require flatness and constant torsion, and not structures corresponding to constant curvature and zero torsion such as the Levi-Civita connection on a Riemannian symmetric space. The key is that the extension of a $\tr$ onto $U(\g)$ allows for a nice algebraic representation of parallel transport, requiring  flatness of the
 connection $\tr$. Indeed, the basic assumption is that $\tr$ extends to the enveloping algebra such that $x\tr(y\tr z)= (x\ast y)\tr z$. From this follows
 \[
 	\llbracket x,y\rrbracket\tr z = (x\ast y-y\ast x)\tr z = x\tr(y\tr z)-y\tr (x\tr z),
\]
and hence $R(x,y,z) = 0$. For \emph{any} connection $\tr$, the corresponding parallel transport of $g$ is 
$$
	g + t f\tr g + \frac{t^2}{2}f\tr(f\tr g) +  \frac{t^3}{3!}f\tr(f\tr(f\tr g)) +\cdots.
$$ 
If the basic assumption above holds, this reduces to the formula $\exp^*(tf)\tr g$.
\end{remark}

Recall that the free Post-Lie algebra over a single generator is the post-Lie algebra of planar rooted trees $\postLie(\{\scalebox{0.6}{\ab}\})$ given in Example~\ref{exm:rootedtrees}. Freeness essentially means that it is a universal model for any post-Lie algebra generated by a single element, and in particular the post-Lie algebra generated by the infinitesimal action of a function $f\in\g^M$ on $M$. For instance, if we decide that $\scalebox{0.6}{\ab}$ represents the element $f \in \g^M$, then there is a unique post-Lie morphism $\F \colon \postLie(\{\scalebox{0.6}{\ab}\})\rightarrow \g^M$ such that $\F(\scalebox{0.6}{\ab}) = f$. Moreover, we then have, e.g.,
$$
	\F(\scalebox{0.6}{\aabb})=\F(\scalebox{0.6}{\ab}\tr\scalebox{0.6}{\ab}) 
	= \F(\scalebox{0.6}{\ab})\tr\F(\scalebox{0.6}{\ab}) = f \tr f,
$$ 
or $\F([\scalebox{0.6}{\ab},\scalebox{0.6}{\aabb}])= [f,f\tr f]$, and so on. This $\F$ is called the elementary differential map, associating planar rooted trees and commutators of these with vector fields on $M$. Hence, all concrete computations in $\g^M$ involving the operations $\tr$ and $[\cdot,\cdot]$ can be lifted to symbolic computations in the free post-Lie algebra $\postLie(\{\scalebox{0.6}{\ab}\})$. 

Revisiting our basic aim, we require a description of $U(\postLie(\{\scalebox{0.6}{\ab}\}))$, which is given as the linear combination of all ordered forests (OF) over the alphabet of planar rooted trees, including the empty forest $\one$, 
\[
	\OF=\Big\{\begin{array}{c}
	\one,\scalebox{0.6}{\ab},\scalebox{0.6}{\ab\ab},\scalebox{0.6}{\aabb},
	\scalebox{0.6}{\ab\ab\ab},\scalebox{0.6}{\aabb\ab},\scalebox{0.6}{\ab\aabb},
	\scalebox{0.6}{\aababb},\scalebox{0.6}{\aaabbb},\scalebox{0.6}{\ab\ab\ab\ab},
	\ldots,\scalebox{0.6}{\aabaabbb},\scalebox{0.6}{\aaabbabb},\ldots
	\end{array}\Big\}.
\]
So an element $a \in U(\postLie(\{\scalebox{0.6}{\ab}\}))$ could be, for instance, of the following form
\[
	a = 3\one + 4.5\scalebox{0.6}{\ab} -2 \scalebox{0.6}{\ab\aabb} + 3 \scalebox{0.6}{\aabb\ab} 
	+ 6 \scalebox{0.6}{\aabaabbb\ab} + 7 \scalebox{0.6}{\ab\aaabbabb} 
	- 2 \scalebox{0.6}{\ab\aaabbabb\ab} \cdots .
\]
To be more precise, $U(\postLie(\{\scalebox{0.6}{\ab}\}))$ consists of all finite linear combinations of this kind, while infinite combinations such as the exponential live in $\widehat{U}\left(\postLie(\{\scalebox{0.6}{\ab}\})\right)$ and are obtained by an inverse limit construction~\cite{FMK}. Elements in the space $\widehat{U}\left(\postLie(\{\scalebox{0.6}{\ab}\})\right)$ we call  Lie--Butcher (LB) series. Note that all computations on such infinite series are done by evaluating the series on something finite in $U(\postLie(\{\scalebox{0.6}{\ab}\}))$. Indeed, we consider $\widehat{U}:=\widehat{U}\left(\postLie(\{\scalebox{0.6}{\ab}\})\right)$ as the (linear) dual space of $U:=U\left(\postLie(\{\scalebox{0.6}{\ab}\})\right)$, with a bilinear pairing $\langle\cdot,\cdot\rangle\colon \widehat{U}\times U\rightarrow \RR$ defined such that $\OF$ is an orthonormal basis, i.e.\ for $\omega,\omega'\in \OF$ we have $\langle \omega,\omega'\rangle = 1$ if $\omega = \omega'$, and zero if $\omega \neq \omega'$.

Two important subclasses of LB-series are
\begin{align*}
	\glb &:= \left\{\alpha\in \widehat{U}(\postLie(\{\ab\}) \ \colon\ \langle \alpha,\one\rangle = 0, 
	\langle\alpha,\omega\shuffle\omega'\rangle = 0\  
	\forall \omega,\omega' \in \OF \backslash\{\one\} \right\}\\
	\Glb &:= \left\{\alpha\in \widehat{U}(\postLie(\{\ab\}) \ \colon\ \langle\alpha,
	\omega\shuffle\omega'\rangle =  \langle\alpha,\omega\rangle 
	\langle\alpha,\omega'\rangle\ \forall \omega,\omega'\in \OF\right\},
\end{align*}
where $\shuffle$ denotes the usual shuffle product of words, e.g., $a \shuffle \one = \one \shuffle a=a$,
$$
	ab \shuffle cd = a(b\shuffle cd) + c(ab \shuffle d).
$$ 
Here elements in $\glb$ are called infinitesimal characters, representing vector fields on $M$ and elements in $\Glb$ are characters, representing flows (diffeomorphisms) on $M$. $\Glb$ forms a group under composition called the \emph{Lie--Butcher group}. A natural question is how does an element $\gamma \in \Glb$ represents a flow on $M$? The elementary differential map sends $\gamma$ to the (formal\footnote{Neglecting convergence of infinite series at this point.}) differential operator, i.e., 
$$
	\F(\gamma) = \sum_{\omega\in\OF} \langle\gamma,\omega\rangle\F(\omega)\in \widehat{U}(\g)^M.
$$ 
The flow map $\Psi_\gamma\colon M\rightarrow M$ is such that the differential operator $\F(\gamma)$ computes the Taylor expansion of a function $\phi \in C^\infty(M,\RR)$ along the flow $\Psi_\gamma$:
\[\F(\gamma)[\phi] = \phi\circ \Psi_\gamma .\]

Recall from Proposition~\ref{prop:Liebracket} that any post-Lie algebra comes with \emph{two} Lie algebras $\g$ and $\overline{\g}$. Hence there are two enveloping algebras $U(\g)$ and $U(\overline{\g})$, with two different associative products. It turns out that $U(\g)$ and $U(\overline{\g})$ are isomorphic as Hopf algebras \cite{EFLMK,EFM}, such that the product of the latter can be represented in $U({\g})$. For LB-series, the resulting two associative products in $U({\g})$ are called the concatenation product and Grossman--Larson (GL) product \cite{GL}. Indeed, we have $\omega \cdot \omega' = \omega\omega'$ (sticking words together). The GL product is somewhat more involved, i.e., for $\alpha,\beta \in \glb$ we have $\alpha\ast \beta = \alpha\cdot\beta + \alpha\tr\beta$, see~\cite{EFLMK} for the general formula.  Interpreted as operations on vector fields on $M$, the GL product represents the standard composition (Lie product) of vector fields as differential operators, while the concatenation represents frozen composition, for $\alpha,\beta\colon M\rightarrow U(\g)$ we have $\alpha\cdot\beta(p) = \alpha(p)\cdot\beta(p)$.  

The two associative products on $U({\g})$ yield two exponential mappings $\exp^\cdot,\exp^\ast$ between $\glb$ and  $\Glb$ obtained from these, 
\[
	\exp^\cdot(\alpha) = \one + \alpha+ \frac12 \alpha\cdot\alpha
		+\frac 16\alpha\cdot\alpha\cdot\alpha+\cdots, 
	\; 
	\exp^\ast(\alpha) = \one + \alpha+ \frac12 \alpha\ast\alpha
		+\frac 16\alpha\ast\alpha\ast\alpha+\cdots .\]
Both send vector fields on $M$ to flows on $M$. However, it turns out that  the Grossman--Larson exponential $\exp^\ast$ sends a vector field to its exact solution flow, while the concatenation exponential $\exp^{\cdot}$ sends a vector field to the exponential Euler flow,
\begin{equation}
	y_1 = \exp(hf(y_0))\cdot y_0.\label{eq:expeuler}
\end{equation}
\smallskip
All the basic Lie group integration methods can be formulated and analysed directly in the space of LB-series $\widehat{U}\left(\postLie(\{\scalebox{0.6}{\ab}\})\right)$ with its two associative products and the lifted post-Lie operation. The Lie-Euler method which moves in successive timesteps along the exponential Euler flow is one such example. A slightly more intricate example is
\begin{example}[Lie midpoint integrator]
On the manifold $M$ a step of the Lie midpoint rule with time step $h$ for~(\ref{eq:ode}),  is given as
\begin{align*}
K &= hf(\exp(K/2)\cdot y_0)\\
y_1 &= \exp(K)\cdot y_0. 
\end{align*}
In $\widehat{U}\left(\postLie(\{\scalebox{0.6}{\ab}\})\right)$, the same integrator $\scalebox{0.6}{\ab}\mapsto \Phi\colon \glb\rightarrow \Glb$ is given as:
\begin{align*}
K &= \exp^\cdot(K/2)\tr (h\scalebox{0.6}{\ab})\in \glb\\
\Phi &= \exp^\cdot(K)\in \Glb. 
\end{align*}
\end{example}

We conclude by commenting that the two exponentials are related exactly by a map $\chi: {\g} \to {\g}$, called post-Lie Magnus expansion \cite{EFLMMK,EFMMK,EFM}, such that $\exp^\cdot(f)=\exp^\ast(\chi(f))$, $f \in {\g}$. The series $\chi(f)$ corresponds to the backward error analysis related to the Lie--Euler method. In $\widehat{U}\left(\postLie(\{\scalebox{0.6}{\ab}\})\right)$ we find
$$
	\chi(\scalebox{0.6}{\ab})
	= \scalebox{0.6}{\ab} 
	- \frac{1}{2} \scalebox{0.6}{\aabb} 
	+ \frac{1}{12}[\scalebox{0.6}{\aabb},\scalebox{0.6}{\ab}] 
	+ \frac{1}{3} \scalebox{0.6}{\aaabbb}
	+ \frac{1}{12} \scalebox{0.6}{\aababb} 
	- \frac{1}{12} \scalebox{0.6}{\aaabbabb} +\cdots
$$
This should be compared with the expansion $\beta$ on page 184 in \cite{LMK}. Chapoton and Patras studied the equality between these exponentials in the context of the free pre-Lie algebra \cite{CP}.


\section*{Acknowledgments}  
The research on this paper was partially supported by the Norwegian Research Council (project 231632).


\ifx\undefined\bysame
\newcommand{\bysame}{\leavevmode\hbox to3em{\hrulefill}\,}
\fi

\end{document}